\documentclass[a4paper, mathscr, 12 pt]{article}
\usepackage{eucal, amssymb,amsmath,graphicx, amsfonts, latexsym}
\def\bold1{\boldsymbol{1}}
\def\bold0{\boldsymbol{0}}

\numberwithin{equation}{section}

\begin{document}
\title{Epidemic modelling: aspects where stochasticity matters}
\author{ Tom Britton, Stockholm University\thanks{Department of
Mathematics, Stockholm University, SE-106 91 Stockholm, Sweden. {\it
E-mail}: tom.britton@math.su.se} \\
David Lindenstrand, Stockholm University\thanks{Department of
Mathematics, Stockholm University, SE-106 91 Stockholm, Sweden. {\it
E-mail}: davlin@math.su.se}\ \thanks{\emph{To whom correspondence
should be addressed.}}}
\date{\today}
\maketitle

\begin{abstract}

\noindent Epidemic models are always simplifications of real world
epidemics. Which real world features to include, and which
simplifications to make, depend both on the disease of interest and
on the purpose of the modelling. In the present paper we discuss
some such purposes for which a \emph{stochastic} model is preferable to a
\emph{deterministic} counterpart. The two main examples illustrate the
importance of allowing the infectious and latent periods to be
random when focus lies on the \emph{probability} of a large
epidemic outbreak and/or on the initial \emph{speed}, or growth
rate, of the epidemic. A consequence of the latter is that estimation
of the basic reproduction number $R_0$ is
sensitive to assumptions about the distributions of the infectious and
latent periods when using the data from the early stages of an
outbreak, which we illustrate with data from the SARS outbreak. Some
further examples are also discussed as 
are some practical consequences related to these stochastic aspects.
\vskip1cm
\end{abstract}

Keywords: stochastic epidemic model, major outbreak probability,
infectious period, latency period, exponential growth rate.

\section{Introduction}

Mathematical epidemic models describe the spread of an
infectious disease in a community (e.g.\ Bailey, 1975, Anderson and
May, 1991, Diekmann and Heesterbeek, 2000). A model can be used to
derive various properties of an outbreak, such as: whether or not a big
outbreak may occur, how big the outbreak will be, or the endemic
level in case the disease becomes endemic. From a
statistical/epidemiological point of view the model and its analysis
may be used to estimate important epidemiological parameters from
observed outbreak data. These estimates can then be used to study
effects of potential interventions to stop or reduce the spreading
of the disease. For example, an endemic disease may go extinct if a
vaccination program is launched having high enough vaccination
coverage (e.g.\ Anderson and May, 1991, pp87, and Gay,
2004), or an outbreak
may be stopped during the early stages of an outbreak if spreading
parameters are reduced enough by means of different sorts of
intervention (e.g.\ Anderson et al., 2004, for an application to SARS).

Mathematical models are always simplifications of reality, but the
hope is that the simplifications have little effect on the epidemic
properties of interest. Simple models have the advantage of being
tractable to analysis and quite often allow for explicit solutions
admitting general qualitative statements. Their main disadvantage is
of course that they may be too simplistic for the conclusions to be
valid also for real world epidemics. Adding more complexity to the
model increases realism but usually makes it harder to analyse and
also introduces more uncertainty by having more parameters. More
complex models are usually analysed by means of numerical solutions
to differential equations, or from numerous stochastic simulations.

The most important features to include to make an epidemic model more
realistic (and at the same time harder to analyse) are to incorporate
individual heterogeneity (e.g.\
Anderson and May, 1991, pp 175) and/or structured mixing patterns (e.g.\
 House and Keeling, 2008, for a deterministic
household model). Another step in
making a model more realistic is to make certain features random,
for example the actual transmission/contact process but also possibly
susceptibility, social structures, the latent period and/or the
infectious period. Such stochastic models thus allow individuals to
behave different from each other in a way that is specified by
random distributions (e.g. Bailey, 1975, Andersson and Britton,
2000a).

Which complexities to include in the model, and which not to,
depend both on the type of disease in question and on the
scientific question motivating the study. The aim of the present
paper is to illustrate some aspects where \emph{stochasticity}
matters. More precisely we focus on two features, the risk for an
outbreak and the initial growth rate of the epidemic, and we
illustrate that they depend heavily on assumptions about the latent
and infectious periods; not only on their mean durations but also on
their randomness. As a consequence, the (stochastic) distribution of
these periods are important when addressing questions relating to
these two features -- using an over-simplified stochastic model or a
deterministic model will give misleading results. For example,
estimating $R_0$ from the initial phase of an epidemic is hard without
additional knowledge about the distributions of the infectious and latents periods, a fact
which we illustrate using data from the SARS outbreak. We illustrate our
results using a simple epidemic model, but the qualitative
conclusions hold also for more realistic models allowing other
heterogeneities. We note that other features of the model, e.g.\ the
basic reproduction number and the outbreak size in case of a major
outbreak, hardly depend on the randomness of
the latent and infectious periods at all, so having a deterministic
latent and infectious period may be appropriate when addressing
other questions.

Most results presented in this paper are not new but have appeared
elsewhere or are "folklore" among stochastic epidemic modellers, but
are perhaps less known outside this community. The aim of the paper
is hence to gather and present the results in a simple form reaching
outside the community of stochastic epidemic modellers. The rest of
the paper is outlined as follows. In Section \ref{model} we present
the standard stochastic SEIR epidemic model for a homogeneously
mixing community of homogeneous individuals. In Section
\ref{properties} properties of the model are presented and
illustrated. In Section \ref{relevance} we interpret the results in
more epidemiologically relevant formulations and illustrate where it
can make a difference. In the discussion we briefly describe, and
give references to, some other situations where stochasticity of
some form affect certain features of the epidemic model.

\section{A simple stochastic epidemic model}\label{model}

\subsection{Definition}

We now define what we call the standard
susceptible-exposed-infectious-removed (SEIR) epidemic model.
Consider a homogeneously mixing community consisting of $n$
homogeneous individuals, where $n$ is assumed to be large. A
transmittable disease is spread according to the following rules.
Initially a small number, $k$, individuals are infectious and the
rest of the community are susceptible to the disease (immune
individuals are simply neglected). Each individual who gets infected
is at first latent (exposed but not yet infectious) for a random
period $L$ with distribution $F_L$. After the latent period has
ended the infectious period starts and lasts for a period $I$ having
distribution $F_I$. All infectious periods and latent periods are
assumed to be mutually independent. While infectious an individual
has random "infectious contacts" at rate $\lambda$, each contact is
with a randomly chosen individual, so the contact rate with a
specific individual is $\lambda/n$ (or more correctly
$\lambda/(n-1)$ but when $n$ is large this distinction is
irrelevant). Contacts with susceptible individuals result in infection
(and their latent period starts); contacts with
non-susceptibles have no effect. Once the infectious period is over
the individual is said to be removed, meaning that the person has
recovered and become immune, and plays no further role in the
epidemic. The epidemic goes on until there are no more infectious or
latent individuals, then the epidemic stops. Let $T$ denote the
(random) number of individuals who get infected during the outbreak,
and that hence are removed at the end of the epidemic.  $T$ is often
called the final size of the epidemic, and $\tilde \rho=T/n$ denotes
the final proportion infected during the outbreak.

In what follows we will restrict ourselves to the case where $L$ and
$I$ have different and independent Gamma distributions, this being a
rather flexible family of distributions. We parametrise these
distributions by their means, $\mu_L=E(L)$ and $\mu_I=E(I)$ ($\geq
0$), and their coefficients of variation $\tau_L=\sqrt{V(L)}/E(L)$
and $\tau_I=\sqrt{V(I)}/E(I)$ $\geq 0$, where $V(\cdot)$ denotes the variance.

\subsection{The basic reproduction number $R_0$}

The perhaps most important property of an epidemic model is the basic
reproduction number, denoted $R_0$, which for the present model can be defined as  the
average number of infections caused by a typical infective when the
disease is introduced into the population. For the present model it
is easy to show that
\[
R_0=\lambda E(I)=\lambda \mu_I.
\]
The basic reproduction number determines both if a major outbreak is
possible, and if so, also the final proportion infected in case there
is a major outbreak. More precisely, it can be shown that $\tilde
\rho$, the ultimate proportion infected, will in a large community
be close to $\rho$, which solves
\begin{equation}\label{eq:rho}
 1-\rho=e^{-R_0 \rho}.
\end{equation}
It is easy to see that $\rho =0$ (corresponding to a minor outbreak)
is always a solution to (\ref{eq:rho}). If $R_0\le 1$, this is in
fact the only solution, meaning that a major outbreak is impossible.
If $R_0>1 $ there is also a unique strictly positive solution
$\rho^*$ ($0< \rho^* < 1$) corresponding to a major outbreak.

As was seen above, $R_0$ only depends on the \emph{mean} of the
infectious period -- not on its randomness nor on the latency
period. The model can be extended to allow for a (perhaps random)
time-varying infectivity $\lambda (s)$ over the infectious period
($0\le s\le I\le \infty$). Then $R_0=E(\int_0^I\lambda (s)ds)$, the
expected accumulated infectivity. As before, $R_0$ determines both
if a major outbreak is possible, and if so, how big the outbreak
will be. In fact, the complete (random) distribution of the final
size, for any finite $n$, can be shown to depend only on the
distribution of the accumulated infectivity $\int_0^I\lambda
(s)ds$ (Ball, 1986), how the infectivity is distributed over time
only affects the time dynamics of the epidemic and not the final
size.

\section{Model properties affected by randomness}\label{properties}

In the previous section it was shown that $R_0$ only depends on the
mean length of the infectious period and not at all on the latent
period. In the present section we study two features, the
probability of a major outbreak and the initial growth rate of the
epidemic, where the randomness of the infectious period and also the latent
period do matter. In the discussion we briefly mention some other
aspects where stochasticity matters.

\subsection{The probability of a major outbreak}\label{prob-outbreak}
When the community $n$ is large, the initial phase of the epidemic
may be approximated by a branching process (Ball, 1986). The reason
for this is that new contacts will most likely be with not yet
contacted people, so new infectives infect (=''give birth'' in
branching process terminology) independently which is the crucial
underlying assumption in branching processes. The branching process
corresponding to our model is the Sevastyanov model
(Jagers, 1975, p 8). Infections correspond to births in the branching
process, the latency period to infancy in the branching process and
the infectious period to the reproductive life stage (life stages
after the reproductive stage play no role for population growth just
like with removed individuals in the epidemic).

Let $\pi$ denote the
probability of a large outbreak (corresponding to
infinite growth of the approximating branching process) when starting
with $k=1$ infectious individual. From
branching process theory it can be shown that $\pi$ is the largest solution to
the balance equation
\begin{equation}\label{eq:pi1}
1- \pi=E[(1-\pi)^X],
\end{equation}
 where $X$ is the (random) number of births of a typical
 individual in the branching process. The balance equation is
 obtained by conditioning on the number of births of the first individual:
 if the first individual has $X=x$ births during her life, all these
 individuals must avoid causing infinite growth.
They do this independently, so the probability for this to happen is $(1-\pi)^x$.

 For our model, with constant
 birth/infection rate $\lambda$ during a Gamma distributed
 infectious period with mean $\mu_I$ and coefficient of variation
 $\tau_I$, the distribution of $X$, and hence also $E[(1-\pi)^X]$,
 can be computed explicitly. By first
 conditioning on the length of the infectious period $I=y$ it is
 easy to show that $X$ then is Poisson distributed with mean
 $\lambda y$, and removing the conditioning makes $X$ follow a
 negative binomial distribution. Using this it can be shown that
 Equation ($\ref{eq:pi1}$) simplifies to
\begin{equation}\label{eq:res1}
1-\pi=\Big(\frac{1}{1+\pi R_0 \tau_I^2}\Big)^{\tau_I^{-2}},
\end{equation}
(this relation can also be found in
Asikainen, 2006, p 28). If for example $\tau_I=0$, implying that the
length of the infectious period is non-random, $\pi$ is the largest
solution to $1-\pi=e^{-\pi R_0}$ which is obtained by taking limits of
(\ref{eq:res1}) when $\tau_I\to 0$. If $\tau_I=1$, corresponding to
an exponentially distributed infectious period, we have that
$\pi=1-1/R_0$ which is clearly different.

By studying the balance equation (\ref{eq:res1}) it is possible to
see how $\pi$, the probability of a major outbreak,  depends on model
parameters. The first conclusion is not very
surprising: $\pi$ is
increasing in $R_0=\lambda \mu_I$ and hence also in the contact rate
$\lambda$ and in the mean infectious period $\mu_I$. A less
obvious conclusion is that $\pi$ is \emph{decreasing} in $\tau_I$
(the coefficient of variation of the infectious period). In other
words, the more random the length of the infectious period is, the
less likely is a major outbreak. Finally, $\pi$ is independent of
$\mu_L$ and $\tau_L$.

\begin{figure}[!h] \begin{center} \bf
\includegraphics[height=!, width=10cm]{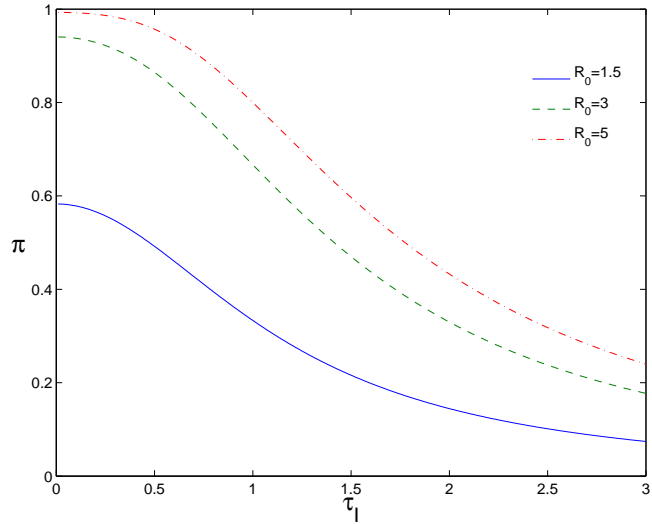}
\caption{\rm The probability of a large outbreak, $\pi$, as function
  of $\tau_I$,
  the coefficient of variation of the infectious period.}
 \label{fig:Pi1}\end{center}
\end{figure}

In Figure \ref{fig:Pi1} we have plotted $\pi$ as a
function of $\tau_I$ in the range 0 (corresponding to a deterministic
infectious period) to 3 (being a very random infectious period), for
three choices of $R_0$.  It is seen that
$\tau_I$ is quite influential. For
example, if $R_0=3$ and $\tau_I=0$, then  $\pi\approx 0.940$. If $R_0$ is
reduced to 1.5 and $\tau_I$ is unchanged we get $\pi\approx 0.583$, whereas
if we instead keep $R_0$ unchanged (at $R_0=3$) and increase
$\tau_I$ to 1, then $\pi\approx 0.667$. It is hence seen that the
variation in the infectious period is as important as $R_0$ for
determining the probability of a major outbreak.

\begin{figure}[!h] \begin{center} \bf
\includegraphics[height=!, width=10cm]{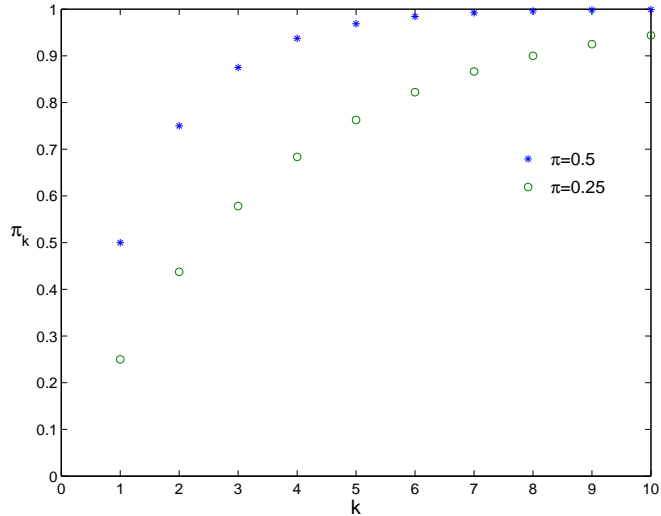}
\caption{\rm The probability $\pi_k$ of an outbreak when $k$ infectious
individuals enters the population.}
\label{fig:Pi2}\end{center}
\end{figure}

The probability $\pi$ defined above was for the case that the epidemic
starts with 1 initially infectious. More generally, we can define
$\pi_k$ as the probability of a major outbreak starting with $k$
initially infectious individuals (so $\pi_1=\pi$). Since, for an
epidemic \emph{not} to take off, none of the initially infectives must
initiate a major outbreak. As a consequence, $\pi_k$ can be expressed
in terms of $\pi_1=\pi$ as
\begin{equation}\label{eq:pm}
\pi_{k}=1-(1-\pi)^{k},
\end{equation}
where $\pi$ is the solution to (\ref{eq:res1}). In Figure
\ref{fig:Pi2}  $\pi_k$ is plotted as a function of $k$ for the cases
$\pi=0.25$ and $\pi=0.5$. It is seen that $\pi_k$ grows quickly up
towards 1, implying that the outbreak probability is close to 1 when
initiated by many individuals as long as $R_0>1$ and $\tau_I$ is not
very large (meaning that infectious
period is not extremely varying).

The distribution of the infectious period is hence mainly of interest
when the epidemic is initiated by rather few individuals.

\subsection{The initial growth rate of the epidemic}\label{sec:growth}
We now study another property which is heavily influenced by both the latent and
infectious periods, their
mean durations as well as their randomness: the initial growth rate of the epidemic. As before
we assume that the community size $n$ is large.

In Section \ref{model} it was shown that an epidemic can only take
off if $R_0=\lambda \mu_I>1$. Since we now focus on the growth
rate of the epidemic we assume this to be the case. As mentioned
before, the early stages of the epidemic in a large community can be
approximated by a branching process. Since $R_0>1$ the branching
process is said to be super-critical, and if the epidemic/branching
process takes off branching process theory (e.g.\ Jagers, 1975)
tells us that the epidemic will grow at an exponential rate during
the initial phase. More precisely, in case of a major outbreak, the number of infectious
individuals at t, $I(t)$, will satisfy $I(t)\sim e^{\alpha t}$
for some $\alpha$. The parameter $\alpha$, denoted the Malthusian
parameter, is known to solve
\begin{equation}\label{eq:Jagers2}
\int_{0}^{\infty}\!\!\!e^{-\alpha t}\lambda P(L<t<L+I)dt=1,
\end{equation}
where $L$ and $I$ are the (random) durations of the latent and
infectious periods respectively. In the present paper $L$ and $I$
are assumed to be independent Gamma distributions with means $\mu_I$
and $\mu_L$, and coefficients of variation $\tau_I$ and $\tau_L$
respectively. The solution $\alpha$ to (\ref{eq:Jagers2}) can
then be shown to solve
\begin{equation}\label{alpha_balance}
 \alpha=\frac{R_0}{\mu_I}\frac{1}{(1+\alpha\tau_L^{2}\mu_L)^{1/\tau_L^{2}}}\left( 1-\frac{1}{(1+\alpha
 \tau_I^{2}\mu_I)^{1/\tau_I^{2}}}\right)
\end{equation}
(see the Appendix for details).

It can be shown that the exponential growth rate $\alpha$ (i.e.\ the
solution to (\ref{alpha_balance})) depends monotonically on all four
parameters of the latent and infectious periods, $\mu_I$, $\mu_L$,
$\tau_I$ and $\tau_L$, keeping $R_0$ fixed. As for the mean
infectious and latent periods, $\mu_I$ and $\mu_L$, the growth is
\emph{decreasing}. This is not surprising: the longer the infectious
period (keeping $R_0=\lambda\mu_I$ fixed!) the slower the epidemic
will grow, and the same applies to the situation where a latent
period becomes longer on average. Perhaps more surprising is that
$\alpha$ depends monotonically on the coefficients of variation
$\tau_L$ and $\tau_I$, and in different ways! It can be shown from
(\ref{alpha_balance}) that the growth rate is \emph{increasing} in
$\tau_L$ but \emph{decreasing} in $\tau_I$. In other words, a more
random latent period increases the growth rate whereas a more random
infectious period decreases it.

A heuristic motivation for the different monotone dependence of the
coefficients of variation goes as follows. Consider first two
alternatives for the infectious period assuming, for simplicity that
there is no latency period: two infectious periods both being two
time-units long (corresponding to small $\tau_I$) and the other
scenario having one infectious period of length 1 and the other of
length 3, thus having the same mean $\mu_I$ but larger $\tau_I$.
During the first time-unit both scenarios will have two persons
infecting but during the second time-unit the first scenario (small
$\tau_I$) will still have two persons infecting, but the second
scenario only one person. During the third time-unit the second
scenario will "catch up" in infecting new people by having one
person infecting (as opposed to no one for the second scenario) but the first scenario will
clearly infect new individuals at an \emph{earlier state in time} thus
resulting in higher growth rate $\alpha$. This motivates why $\alpha$
is decreases in $\tau_I$. The motivation for the
growth rate being increasing in $\tau_L$ is similar. Suppose that we have
two alternative scenarios similar to before: two latent periods of
equal length two time-units, or one of length 1 and one of
length 3, and assume for simplicity that all infectious periods last one time-unit in
both scenarios. In the first scenario the two individuals will
infect others between time 2 and 3 whereas in the second scenario
one person will infect between time 1 and 2 and the other between 3
and 4. The second scenario (with higher $\tau_L$) will have a higher
growth rate because of the multiplicative effect the first person's
infections will cause: these people will start new epidemic outbreaks at
an earlier state.

\begin{figure}[!h] \begin{center} \bf
\includegraphics[height=!, width=14cm]{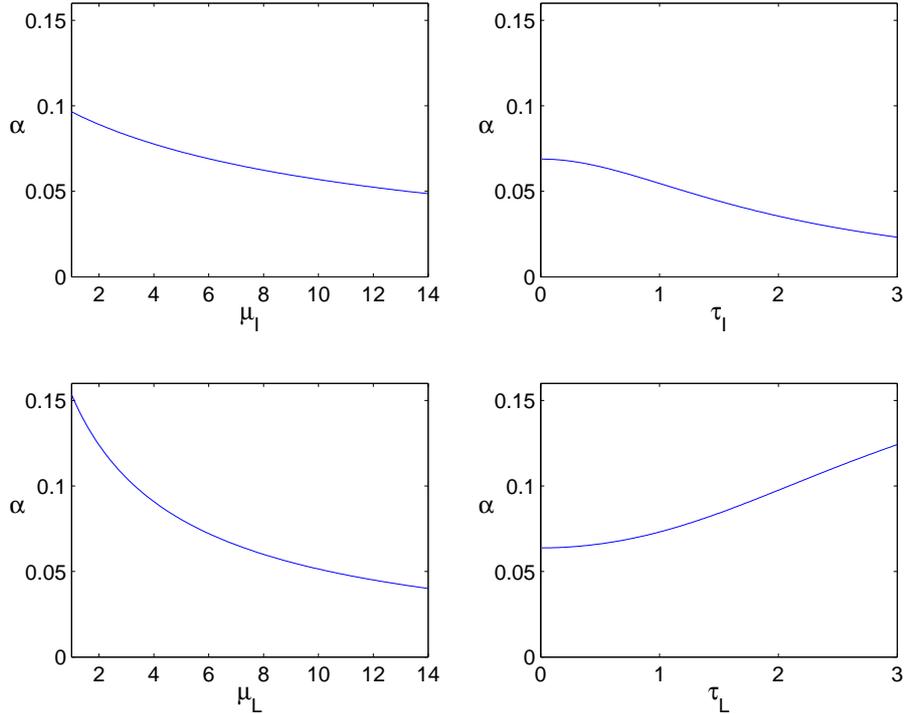}
\caption{\rm Plots of the initial exponential (per day) growth rate $\alpha$ as
  function of the model parameters. Parameters not varied are set to:
  $R_0=2$, $\mu_L=\mu_I=7$ $\tau_L=\tau_I=\frac{3}{7}$.}
\label{fig:alpha2}\end{center}
\end{figure}

In Figure \ref{fig:alpha2} the exponential growth rate $\alpha$ is
computed numerically for the case $R_0=2$, this being a common value
for diseases like influenza (e.g.\ Mills et al., 2004). In each of the four sub-plots, one
parameter is varied over an interval (1 to 14 days for the mean
durations $\mu_L$ and $\mu_I$ and 0 to 3 for the coefficients of
variation $\tau_L$ and $\tau_I$) keeping the remaining parameters
constant. The means are set to 7 days and the coefficients of
variation to 3/7 (corresponding to a standard deviation of 3 days)
when not varied.

From the figure it is clear that all four parameters $\mu_L$,
$\mu_I$, $\tau_L$ and $\tau_I$ are quite influential for the initial
growth rate of the epidemic. As mentioned above, the growth rate is
decreasing in the two mean durations (recall that the expected
accumulated infectivity $R_0=\lambda\mu_I$ is kept fixed, so when
$\mu_I$ changes, so does $\lambda$). As for the coefficients of
variation, the growth rate $\alpha$ decreases with $\tau_I$ but
increases with $\tau_L$.

\section{Practical relevance}\label{relevance}

\subsection{Estimating $R_0$ from growth rate needs prior knowledge}
Our first, and perhaps most important observation, lies in the
consequences of knowing that the growth rate depends heavily on all
of the parameters $\mu_L$, $\mu_I$, $\tau_L$ and $\tau_I$, and not
only $R_0$. This implies that it is harder to estimate $R_0$ from
only observing the early stages of an epidemic as we now illustrate.

Recently, an important area in infectious disease epidemiology has
been to analyse emerging infectious diseases, for example SARS
(e.g.\ McLean et al., 2005) and the fear for a pandemic influenza
(e.g.\ Ferguson et al., 2004). One important task when analysing
emerging infectious diseases is to estimate $R_0$ using data from
the initial phase of the epidemic. Such data sets typically consists
of the number of diagnosed cases (per day or per week) over a
certain observation period, typically weeks or months. One can argue
that the number of diagnosed cases roughly corresponds to the number
of recovered individuals, and using branching process theory it can
be shown that this number will have the same growth rate as the
number of infectives. If we let $R(t)$ denote the accumulated number
of removed individuals up to time $t$, it is known from branching
process theory that
\begin{equation}\label{growth}
R(t)\approx W e^{\alpha t},
\end{equation}
where $W$ is random variable, the same for all $t$, and $\alpha$ is
the Malthusian parameter treated in Section 3.2. If we look at the
ratio of the number of removed individuals for two different
observation times it follows that $R(t_1)/R(t_0)\approx e^{\alpha
  (t_1-t_0)}$ implying that we can estimate the growth rate $\alpha$
by
\begin{equation}
\hat \alpha =\frac{\log (R(t_1))-\log (R(t_0))}{t_1-t_0}.
\label{alpha-est}
\end{equation}
The time-points $t_0<t_1$ should be chosen such that the epidemic
has really taken off at $t_0$ and not too many should have been
infected by $t_1$.

The remaining problem lies in making
conclusions about $R_0$ from the estimate $\hat \alpha$. Equation
(\ref{alpha_balance}) gives a one-to-one correspondence between $\alpha$
and $R_0$ when the model parameters for the latent and infectious
period are given. Rearranging Equation (\ref{alpha_balance}) gives the
following expression for $R_0$:
\begin{equation}
R_0=\alpha
\mu_I\frac{(1+\alpha\tau_L^{2}\mu_L)^{1/\tau_L^{2}}}{\left(1-\frac{1}{(1+\alpha
 \tau_I^{2}\mu_I)^{1/\tau_I^{2}}}\right)}
\end{equation}
However, for emerging infectious diseases the parameters of the
infectious and latent periods are rarely known. The best one can
hope for are some crude estimates. This will induce uncertainty
in the estimate for $R_0$ no matter how precise the estimator $\hat
\alpha$ is.

\begin{figure}[!h] \begin{center} \bf
\includegraphics[height=!, width=10cm]{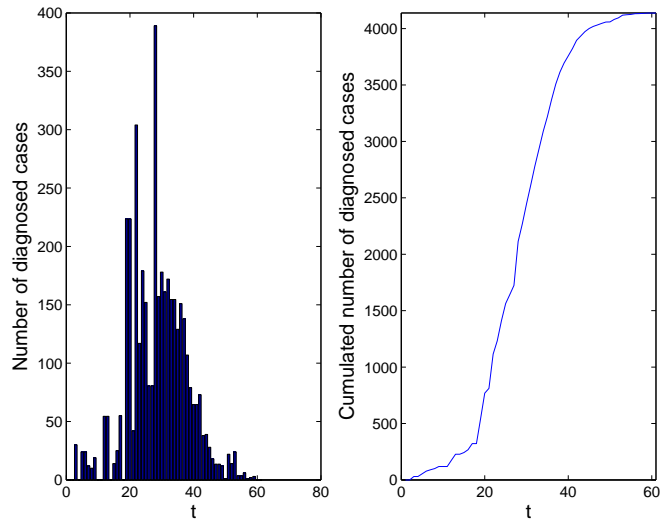}
 \caption{\rm Sars outbreak in China 2003.04.02 - 2003.06.02. Data from WHO.}
 \label{fig:sars}\end{center}
\end{figure}

We now illustrate this using WHO data from the SARS outbreak (WHO
webpage). Our model is of course unrealistic for this outbreak in
several aspects as we are neglecting other community
heterogeneities. However, the same qualitative conclusions would
hold also for more realistic models. In Figure \ref{fig:sars} part
of a large outbreak of SARS in China is illustrated. It shows the
incidence and accumulated number of diagnosed SARS cases by the day,
between April and June in 2003.

From this data we estimate the growth rate $\alpha$ using
(\ref{alpha-est}). Rather than estimating $\alpha$ from one time
interval $(t_0,\ t_1)$ we take several, thus getting several
$\alpha$-estimates. We then take the mean of these estimates as our
final estimate. More precisely we took the intervals $(t_0, t_1)$=
(10,20), (10,25), and (15,25), all three representing the early
stages of the epidemic neglecting the very first bit and stopping
before the speed really starts dropping. The resulting
$\alpha$-estimates were $\hat \alpha_1=0.071$, $\hat \alpha_2=0.054$, and
$\hat \alpha_3=0.034$. We take the mean of these values as our final
estimate: $\hat \alpha=0.0530$. (We will use the estimate to
illustrate that a range $R_0$-values are 
consistent with this estimate, the exact value of $\hat\alpha$ is of
secondary importance.)

Given the estimate ($\hat\alpha=0.0530$) we now use Equation
(\ref{alpha_balance}) to see what we can say about $R_0$. The
disappointing answer is that, unless we assume some prior knowledge
about the latent and infectious periods, we can hardly say anything
about $R_0$, except that $R_0>1$ since the epidemic is taking off.
In order to say more about $R_0$ one needs either more detailed data
or some other knowledge about the latent and infectious periods. If
infections are contact-traced it is possible to make inference on
the generation times. However estimating model parameters from such
inference is far from simple (Svensson, 2007). If such information
is not available, $R_0$ can be estimated by assuming interval ranges
for each model parameter, ranges within which the true parameter
values are believed to lie. To illustrate this from the SARS data we
choose the following intervals: $\mu_I$ and $\mu_L$ is assumed to
lie between 3 and 11 days (with 7 days as mid-point), and the
coefficients of variation are assumed to lie between $0$ and $4/7$
(corresponding to $4$ days for the mid-points above). In Table
$\ref{table:R0}$ the $R_0$ estimate, based on (\ref{alpha_balance}),
$\hat\alpha =0.053$ and current values of $\mu_L$, $\mu_I$, $\tau_L$
and $\tau_I$, is listed for each of the 16 combinations interval
end-points. The point estimate when each parameter takes on the
mid-value ($\mu_L=\mu_I=7$ and $\tau_L=\tau_I=2/7$ and $\hat
\alpha=0.053$) equals $\hat R_0=1.747$.

\begin{table}[!h]
\begin{center}
\setlength\tabcolsep{7pt}
\begin{tabular}{c c c c | c}
   $\mu_L$ & $\mu_I$ & $\tau_L$ & $\tau_I$ & $\hat R_0$\\
\hline
 $3$ & $3$ & $0$ & $0$ & $1.2903$\\
 $3$ & $3$ & $0$ & $4/7$ & $1.2935$\\
 $3$ & $3$ & $4/7$ & $0$ & $1.2897$\\
 $3$ & $3$ & $4/7$ & $4/7$ & $1.2930$\\

 $3$ & $11$ & $0$ & $0$ & $1.5528$\\
 $3$ & $11$ & $0$ & $4/7$ & $1.6468$\\
 $3$ & $11$ & $4/7$ & $0$ & $1.5521$\\
 $3$ & $11$ & $4/7$ & $4/7$ & $1.6461$\\

 $11$ & $3$ & $0$ & $0$ & $1.9674$\\
 $11$ & $3$ & $0$ & $4/7$ & $1.9724$\\
 $11$ & $3$ & $4/7$ & $0$ & $1.8834$\\
 $11$ & $3$ & $4/7$ & $4/7$ & $1.8881$\\

 $11$ & $11$ & $0$ & $0$ & $2.3677$\\
 $11$ & $11$ & $0$ & $4/7$ & $2.5111$\\
 $11$ & $11$ & $4/7$ & $0$ & $2.2666$\\
 $11$ & $11$ & $4/7$ & $4/7$ & $2.4039$\\
\end{tabular}
\end{center}
\caption{Estimates of $R_0$ from SARS outbreak in China 2003.04.02 -
2003.06.02. within different assumptions of model parameters and
$\hat\alpha=0.053$. Data from WHO.}\label{table:R0}
\end{table}

As can be seen from the table the estimate $\hat R_0$ depends quite
a lot on our assumptions about the latent and infectious periods.
The smallest estimate is $\hat R_0=1.2897$ obtained when $\mu_L$,
$\mu_I$ and $\tau_I$ are at their minimal possible point and where
$\tau_L$ is at its maximal point. The largest estimate is
$R_0=2.5111$ obtained for the ``opposite'' parameter choices. Within
the range of "possible" parameter values for the latent and
infectious periods, the $R_0$ estimate hence changes by a factor 2.
It is hence hard to make precise estimates of $R_0$ without other
sources of information regarding the latent and infectious periods.

This illustrates that an estimate of $R_0$ using data from
the initial growth is quite uncertain except in the rare case
that the parameters of the latent and infectious periods are known
with fairly high precision.

\subsection{Estimating variability from final size data}

In Section \ref{prob-outbreak} it was shown that, for fixed $R_0$,
the more random the infectious period is, the more unlikely is a
large outbreak (the same conclusion holds when other factors, e.g.\
susceptibilities and/or infectivities, are varied, Andersson and
Britton, 2000a and references therein). This observation can be used
to say something about the randomness of the infectious period
(and/or of individuals) from final size data, i.e.\ data lacking any
time measurements. If we observe the final proportion infected
$\tilde \rho$ in a large outbreak we estimate $R_0$ using Equation
(\ref{eq:rho}), which gives the estimate
\begin{equation}
\hat R_0=\frac{-\ln (1-\tilde \rho)}{\tilde \rho}.
\end{equation}
The information about $\tau_I$ lies in the fact that a major
outbreak took place, an event with small probability when $\tau_I$
is large. In the Bayesian framework this can be illustrated by
comparing the prior distribution $p(\tau_I)$ with the posterior
distribution $p(\tau_I| \tilde \rho)$. Using Bayes formula we get
\[
p(\tau_I| \tilde \rho) \propto p(\tilde \rho|\tau_I)p(\tau_I),
\]
i.e.\ the posterior distribution equals the prior distribution
multiplied by $p(\tilde \rho|\tau_I)$, the probability of a large
outbreak, denoted $\pi$ in Section \ref{prob-outbreak}. There it was
shown that $\pi=p(\tilde \rho|\tau_I)$ was decreasing in $\tau_I$,
the coefficient of variation of the infectious period. So, any prior
knowledge about $\tau_I$ is shifted towards smaller values in the
posterior distribution. The same type of conclusion also applies to
other individual heterogeneities: the fact that a major outbreak has
occurred shifts any prior knowledge about individual variation
towards less variation.

\begin{figure}[!h] \begin{center} \bf
\includegraphics[height=!, width=10cm]{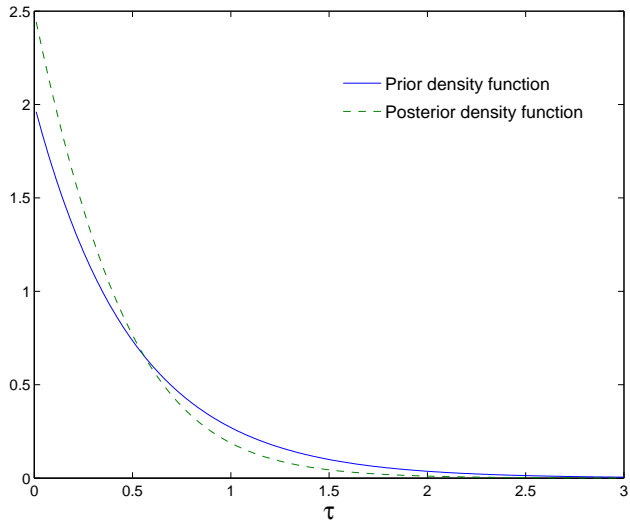}
 \caption{\rm Exponential prior distribution function of $\tau_I$ with
   mean $0.5$, and posterior distribution of $\tau_I$
after observing an outbreak resulting in 50\% getting infected, with mean $0.384$.}
 \label{fig:dist}\end{center}
\end{figure}

Of course, more detailed data containing time-measurements, or data
from more than one outbreak is to be preferred. But, if no such data
is available, any prior information about the randomness of the
infectious period is shifted towards smaller values of $\tau_I$ when
inference is based on final size data from one major outbreak. In
Figure \ref{fig:dist} we illustrate this for the case that $\tau_I$
has an exponential distribution with mean 0.5 as prior distribution,
and where the posterior distribution is based on a major outbreak
resulting in 50\% getting infected.

It is seen that the posterior distribution is not more
"concentrated" compared with the prior distribution, as is usually
the case. Instead the posterior distribution is merely shifted
towards smaller values implying that the belief after observing an
epidemic with 50\% getting infected results in that the posterior
favors smaller values of $\tau_I$ (the coefficient of variation of
the infectious period), as compared to prior beliefs. The posterior
mean of $\tau_I$ is ??, The same qualitative conclusion, that a
major outbreak results in higher posterior belief for small
coefficient of variation of the infectious period, holds for any
fraction getting infected and any prior distribution for $\tau_I$.

\section{Discussion}

In the present paper we have tried to motivate the use of stochastic
models when studying certain features in epidemics. First it was
illustrated that the \emph{probability} for a major outbreak is
greatly affected by the randomness of the infectious period, or more
generally, the randomness of the ``infectivity'' exerted by an
individual.  The more variation the distribution of the infectious
period contains, the less likely is a major outbreak. As a
consequence, observed epidemics (major outbreaks) will tend to
originate from diseases with infectious periods not having very
skew/heavy-tailed distributions. It was also shown that the
probability of a major outbreak is unaffected by a latency period of
arbitrary length. The latter result relies on the assumption that
individuals do not change behaviour as the epidemic progresses nor
that preventive measures are put into place -- then a latency period
will have an effect.

The second feature studied was the initial exponential growth rate.
This rate was shown to depend heavily on both the latent and
infectious periods, there means as well as their randomness. From a
practical perspective
 this implies that, unless additional
information about the infectious period and latency period
distributions is available, it is very hard to estimate
the basic reproduction number $R_0$ (and effects of possible preventive
measures) from the exponential growth rate of the
initial outbreak phase.

There are also other features in epidemics affected by randomness
and not only mean values. Common for most of these situations are
that, for some type of event, only few random objects are
influential. One such feature is the time to disease extinction of
endemic diseases: before disease extinction only few are infectious.
For example, Andersson and Britton (2000b) show that not only the
means but also the coefficients of variation of the
 latency period, infectious period and life-duration affect the time to
extinction when starting at the endemic level.

Another feature affected by randomness is \emph{vaccine response}.
Two models for vaccine response are the \emph{leaky} model and the
\emph{all-or-nothing} model (Halloran et al., 1992). The leaky model
assumes that each person vaccinated has a susceptibility that is
reduced by a factor $e$ (for efficacy). The all-or-nothing model
instead assumes that a proportion $e$ are completely immune whereas
the remaining proportion vaccinated are unaffected by the vaccine.
Here too, $e$ is called efficacy. In both cases, the relative risk
that a vaccinated person gets infected by an infectious contact is
$1-e$ (so the person avoids infection due to the vaccine with
probability $e$). Even though the two models have the same
''efficacy'' their effect is different. In fact, a leaky vaccine
always reduces the spread less than an all-or-nothing vaccine with
the same efficacy -- so the randomness in vaccine effect matters. A
simple explanation to this is the following (Ball and Becker, 2006).
Both vaccine models have the same probability of infection ($1-e$)
at the first contact with an individual. However, among the
vaccinated people who escape infection upon the first contact,
people vaccinated with a leaky vaccine still have relative
susceptibility $1-e$ whereas those with the all-or-nothing vaccine
escaping infection the first time all have the ''all''-effect and
are hence completely immune. As a consequence, the final size in
case of a major outbreak will be smaller with an all-or-nothing
vaccine as compared to a leaky vaccine having the same efficacy $e$
($0\le e\le 1$). Still, both vaccine responses have the same
critical vaccination coverage $v_c=e^{-1}\left(1-1/R_0\right)$,
meaning that the same fraction has to be vaccinated with either
vaccine in order to obtain herd immunity.

As pointed out there are many other features less influenced by
stochasticity, for example $R_0$. In the present paper we simply
focus on aspects where stochasticity does matter.

Needless to say, the model we have studied is by no means fully
realistic. Important extensions are for example to allow for
different types of individuals having different susceptibility,
infectivity and/or mixing patterns, e.g.\ households with higher
contact rates within households (since households are small,
stochasticity play a roll also here, cf.\ Ball et al., 1997).
However, the features considered in the present paper are still
valid under such more realistic models.

\appendix
\section*{Appendix: The Malthusian parameter}\label{sec:app2} The Malthusian
parameter $\alpha$ is the solution to (\ref{eq:Jagers2}). To begin
with,
\begin{eqnarray*}
\lambda P(L<t<L+I) &=&\lambda \int_{0}^{t} f_L(s)\int_{t-s}^{\infty}f_I(r)dr ds\\
    &=& \lambda \int_{0}^{t} f_L(s)(1-F_I(t-s))ds.
\end{eqnarray*}
Hence (\ref{eq:Jagers2}) equals
\begin{eqnarray}
 \lambda \int_{0}^{\infty}e^{-\alpha t} \int_{0}^{t} f_L(s)(1-F_I(t-s))ds dt
        &=& \lambda \int_{0}^{\infty}e^{-\alpha s} f_L(s)
         \int_{s}^{\infty}e^{-\alpha (t-s)} (1-F_I(t-s)) dt ds \nonumber\\
        &=& \lambda \varphi_L(\alpha)(1-\varphi_I(\alpha))\frac{1}{\alpha}.\label{eq:jagers3}
\end{eqnarray}
The second equality follows from partial integration and identifying
the laplace transforms of the latent and infectious periods:
$\varphi_L(\alpha)=E(e^{-\alpha L})=\int_0^\infty e^{\alpha
s}f_L(s)ds$, and similar for the infectious period. The infectious
period $I$ is gamma-distributed. Using first the more common
parametrization $I\sim \Gamma (\alpha_I,\beta_I)$ we get
\begin{equation}\label{eq:jagers4}
 \varphi_I(\alpha)=\Big(\frac{\beta_I}{\beta_I+\alpha}\Big)^{\alpha_I}.
\end{equation}
Since also the latent period is gamma distributed we also have that
$\varphi_L(\alpha)=\Big(\frac{\beta_L}{\beta_L+\alpha}\Big)^{\alpha_L}$.
Thus, by $(\ref{eq:jagers3})$ and $(\ref{eq:jagers4})$, Equation
(\ref{eq:Jagers2}) is simplified to
\begin{equation}\label{app2}
 \alpha=\lambda \Big(\frac{\beta_L}{\beta_L+\alpha}\Big)^{\alpha_L} \Big(1-\Big(\frac{\beta_I}{\beta_I+\alpha}\Big)^{\alpha_I}\Big).
\end{equation}
If we convert to the more interpretable parameters mean $\mu$ and
coefficient of variation $\tau$ we have that
$\mu_I=\alpha_I/\beta_I$ and $\tau_I=1/\sqrt{\alpha_I}$, and
similarly for the latent period. Equation (\ref{app2}) can then,
after some simple algebra and using that $R_0=\lambda \mu_I$, be
written as
\begin{equation*}\label{app3}
 \alpha=\frac{R_0}{\mu_I}\frac{1}{(1+\alpha\tau_L^{2}\mu_L)^{1/\tau_L^{2}}}\left( 1-\frac{1}{(1+\alpha
 \tau_I^{2}\mu_I)^{1/\tau_I^{2}}}\right).
\end{equation*}

\section*{Acknowledgements}

We thank {\AA}ke Svensson for help in simplifying formulae for the Malthusian parameter.
T.B. is grateful to the Swedish Research Council for financial
support.


\begin{thebibliography}{99}

\bibitem{af04} Anderson, R.M., Fraser, C., Ghani, A.C., Donnelly,
  C.A., Riley, S., Ferguson, N.M., Leung, G.M., Lam, T.H., and Hedley,
  A.J.\
  (2004). Epidemiology, transmission dynamics and control of SARS: the
  2002-2003 epidemic. {\it Philos Trans R Soc Lond B Biol
  Sci.\bf 359}, 1091-1105.

\bibitem{am91} Anderson R.M.\ and May R.M.\ (1991). {\it Infectious diseases of
humans; dynamic and control}. Oxford: Oxford University Press.

\bibitem{ab00a} Andersson, H.\ and Britton, T. (2000a): Stochastic epidemic models
and their statistical analysis. {\it Springer Lecture Notes in
Statistics, \bf 151}. Springer-Verlag, New York.

\bibitem{ab00b} Andersson, H.\ and Britton, T.\ (2000b): Stochastic epidemics in
dynamic populations: quasi-stationarity and extinction.
{\it J. Math. Biol., \bf 41}, 559-580.

\bibitem{a06} Asikainen, T.\ (2006). Some results in the field
of epidemic modeling and analysis of a smallpox outbreak. Research
Report 2006:5, Mathematical Statistics, Stockholm University.

\bibitem{b75} Bailey,  N.T.J.\ (1975). {\it The Mathematical Theory of
Infectious Diseases and its Applications}. London: Griffin.

\bibitem{b86} Ball, F.G.\ (1986). A unified approach to the distribution
of total size and total area under the trajectory of the infectives
in epidemic models. {\it Adv. Appl. Prob. \bf 18}, 289-310.

\bibitem{BB06} Ball, F.G.\ and Becker, N.G.\ (2006). Control of transmission with two
types of infection. {\em Math. Biosci.\/} {\bf 200}, 170-187.

\bibitem{bms97} Ball F., Mollison D. and Scalia Tomba G. (1997): Epidemics in populations with two levels
of mixing {\it Ann. Appl. Prob. \bf 7}, 46-89.

\bibitem{DH00} Diekmann, O.\  and Heesterbeek, J.A.P.\ (2000).
{\it Mathematical epidemiology of infectious diseases: model
building, analysis and interpretation}. Chichester: John Wiley

\bibitem{nige4} Gay N.J.\ (2004). The theory of measles
  elimination: implications for the design of elimination
  strategies. {\it The Journal of infectious diseases, \bf 189},
  S27-S35.

\bibitem{F06} Ferguson, N.M., Cummings, D.A., Fraser, C., Cajka, J.C.,
Cooley, P.C., Burke, D.S. (2006). Strategies for mitigating an
influenza pandemic. {\it Nature,\bf 442}, 448-452.

\bibitem{HHL92} Halloran, M.E., Haber, M., Longini, I.M.\ (1992)
Interpretation and estimation of vaccine efficacy under
heterogeneity. {\it
Amer.~J.~Epidemiol., \bf 136}, 328-343.

\bibitem{House08} House T. and Keeling M.J. (2008). Deterministic epidemic
models with explicit household structure.{\it Math. Biosci.  \bf 213}, 29-39.

\bibitem{test}Jagers P. (1975). {\em Branching Processes with
Biological Applications}. John Wiley: London.

\bibitem{mcle05} McLean, A., May, R., Pattison, J., Weiss, R. (eds.) (2005).
{\it SARS: a case study in emerging infections}
Oxford: Oxford University Press.

\bibitem{Rob04} Mills, C.E., Robins, J.M., and Lipsitch, M.\ (2004).
Transmissibility of 1918 pandemic influenza. {\it Nature \bf 432},
904-906,


\bibitem{S07} Svensson, {\AA}. (2007). A note on generation times
in epidemic models. {\it Math. Biosci. \bf 208}, 300-311.

\bibitem{who} WHO. Cumulative Number of Reported Probable Cases of Severe Acute Respiratory Syndrome (SARS).
Available from: URL: http://www.who.int/csr/sars/country/en/.

\end{thebibliography}
\end{document}